# FLEXIBLE SUSPENSIONS WITH A HEXAGONAL EQUATOR

VICTOR ALEXANDROV AND ROBERT CONNELLY

Sobolev Institute of Mathematics, Novosibirsk, Russia,
Novosibirsk State University, Novosibirsk, Russia,
and Cornell University, Ithaca, NY, USA

May 22, 2009

ABSTRACT. We construct a flexible (non immersed) suspension with a hexagonal equator in Euclidean 3-space and study its properties related to the Strong Bellows Conjecture which reads as follows: if an immersed polyhedron $\mathcal{P}$ in Euclidean 3-space is obtained from another immersed polyhedron $\mathcal{Q}$ by a continuous flex then $\mathcal{P}$ and $\mathcal{Q}$ are scissors congruent.

## § 1. INTRODUCTION

A polyhedron (more precisely, a polyhedral surface) is said to be *flexible* if its spatial shape can be changed continuously due to changes of its dihedral angles only, i.e., in such a way that every face remains congruent to itself during the flex.

For the first time flexible sphere-homeomorphic polyhedra in Euclidean 3-space were constructed by R. Connelly in 1976 [8]. Since that time, many properties of flexible polyhedra were discovered, for example:

(1) In 1985 R. Alexander [1] proved that *every flexible polyhedron in Euclidean n-space, $n \geq 3$, preserves the following quantity*

$$\sum \varphi(F) \operatorname{vol}_{n-2}(F)$$

called its *total mean curvature*. Here summation is taken over all $(n-2)$-faces $F$ of the polyhedron, $\operatorname{vol}_{n-2}(F)$ denotes the $(n-2)$-volume of $F$, and $\varphi(F)$ denotes the dihedral angle between the two $(n-1)$-faces adjacent to $F$. Later it was shown by several authors that invariance of the total mean curvature is a consequence of the Schläfly differential formula, see, e. g., [3].

(2) In 1996 I.Kh. Sabitov [15] proved that that *every flexible polyhedron in Euclidean 3-space preserves its oriented volume during the flex* and, thus, gave a negative answer to the Bellows Conjecture. An improved presentation is given in [16]; another proof is published in [9].

1991 *Mathematics Subject Classification.* 52C25.

The first author receives support from the Russian State Program for Leading Scientific Schools via Grant NSh–8526.2008.1. Research of the second author is supported in part by NSF Grant No. DMS–0209595.

Typeset by $\mathcal{A}_{\mathcal{M}}\mathcal{S}$-TEX

In this paper we construct a new example of a flexible polyhedron (with self-intersections) in Euclidean 3-space, namely, a flexible suspension with a hexagonal equator, and study its properties related to the Strong Bellows Conjecture which reads as follows: *if an immersed polyhedron $\mathcal{P}$ in Euclidean 3-space is obtained from another immersed polyhedron $\mathcal{Q}$ by a continuous flex then $\mathcal{P}$ and $\mathcal{Q}$ are scissors congruent*, i.e., $\mathcal{P}$ can be divided into a finite set of tetrahedra each of which can be moved independently one from another in space in such a way that the resulting set will give a partition of $\mathcal{Q}$.

## §2. A description of flexible suspensions

Our study of flexible suspensions is based on the fundamental flexibility equation derived in [7]. In this section we briefly review the notation and facts which are important for us.

**Basic definitions.**

A map from a simplicial complex $K$ to Euclidean 3-space $\mathbb{R}^3$, linear on each simplex of $K$, is called a *polyhedron*. If the vertices of $K$ are $v_1, \ldots, v_V$, and if $\mathcal{P} : K \to \mathbb{R}^3$ is a polyhedron, then $\mathcal{P}$ is determined by the $V$ points $p_1, \ldots, p_V$, called the *vertices* of $\mathcal{P}$, where $\mathcal{P}(v_j) = p_j$.

If $\mathcal{P} : K \to \mathbb{R}^3$ and $\mathcal{Q} : K \to \mathbb{R}^3$ are two polyhedra then we say $\mathcal{P}$ and $\mathcal{Q}$ are *congruent* if there is an isometry $A : \mathbb{R}^3 \to \mathbb{R}^3$ such that $\mathcal{Q} = A \circ \mathcal{P}$, i.e., which takes each vertex of $\mathcal{P}$ to the corresponding vertex of $\mathcal{Q}$, $q_j = A(p_j)$ or equivalently $\mathcal{Q}(v_j) = A(\mathcal{P}(v_j))$ for all $j = 1, \ldots, V$. We say $\mathcal{P}$ and $\mathcal{Q}$ are *isometric* if each edge of $\mathcal{P}$ has the same length as the corresponding edge of $\mathcal{Q}$, i.e., if $\langle v_j, v_k \rangle$ is a 1-simplex of $K$ then $|p_j - p_k| = |q_j - q_k|$, where $|\cdot|$ stands for the Euclidean norm in $\mathbb{R}^3$.

A polyhedron $\mathcal{P}$ is *flexible* if, for some continuous one parameter family of polyhedra, $\mathcal{P}_t$, $0 \leqslant t \leqslant 1$, the following three conditions hold true: (1) $\mathcal{P}_0 = \mathcal{P}$; (2) each $\mathcal{P}_t$ is isometric to $\mathcal{P}_0$; (3) some $\mathcal{P}_t$ is not congruent to $\mathcal{P}_0$.

Let $K$ be defined as follows: $K$ has vertices $v_0, v_1, \ldots, v_n, v_{n+1}$, where $v_1, \ldots, v_n$ form a cycle ($v_j$ adjacent to $v_{j+1}$, $j = 1, \ldots, n-1$ and $v_n$ adjacent to $v_1$) and $v_0$ and $v_{n+1}$ are each adjacent to all of $v_1, \ldots, v_n$. Call $\mathcal{P}(v_0) = N$ the north pole and $\mathcal{P}(v_{n+1}) = S$ the south pole and $\mathcal{P}(v_j) = p_j$, $j = 1, \ldots, n$, vertices of the *equator*. Such a $\mathcal{P}$ is called a *suspension*, see Fig. 1.

**The variables and the fundamental equation.**

Let $\mathcal{P} : K \to \mathbb{R}^3$ be a suspension; $N$ and $S$ be the north and south poles of $\mathcal{P}$; and $p_1, p, \ldots, p_n$ be the vertices on the equator in cyclic order, see Fig. 1. Below we use the sign $\stackrel{\text{def}}{=}$ as an abbreviation of the phrase 'put by definition'.

$$e_j \stackrel{\text{def}}{=} N - p_j, \quad e'_j \stackrel{\text{def}}{=} p_j - S, \quad j = 1, 2, \ldots, n, \quad e_{n+1} \stackrel{\text{def}}{=} e_1, e'_{n+1} \stackrel{\text{def}}{=} e'_1,$$

$$e_{j,j+1} \stackrel{\text{def}}{=} e_j - e_{j+1} = e'_{j+1} - e'_j \quad \text{are edges of the equator,}$$

$$R \stackrel{\text{def}}{=} e_j + e'_j = N - S, \quad x \stackrel{\text{def}}{=} R \cdot R,$$

where $\cdot$ stands for the inner product in $\mathbb{R}^3$.

It is easy to see that if a suspension is flexible then $x$, the squared distance between the north and south poles, is non-constant. In the sequel we treat $x$ as an independent variable and consider all other expressions as functions of $x$.

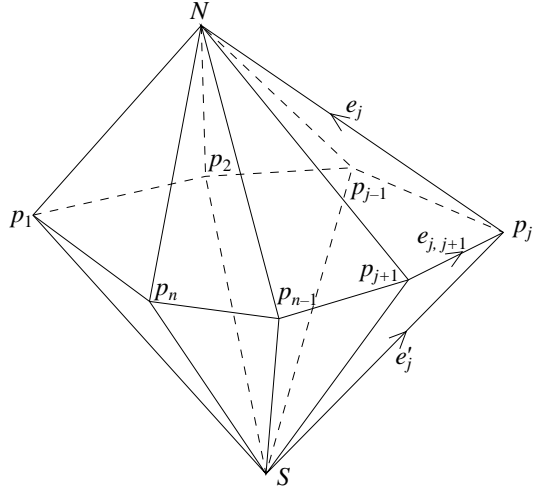

Fig. 1.

Let $\pi : \mathbb{R}^3 \to \pi(\mathbb{R}^3)$ be the orthogonal projection onto the plane perpendicular to $R$ which is regarded as the plane of complex numbers. Using standard facts from analytic geometry and provided $\pi(e_j) \neq 0$ for all $j = 1, 2, \ldots, n$, we express the angle $\theta_{j,j+1}$ from $\pi(e_j)$ to $\pi(e_{j+1})$ in $\pi(\mathbb{R}^3)$ by the formula

$$e^{i\theta_{j,j+1}} = \frac{\pi(e_j)}{|\pi(e_j)|} \cdot \frac{\pi(e_{j+1})}{|\pi(e_{j+1})|} = \frac{G_{j,j+1}}{H_j H_{j+1}} \stackrel{\text{def}}{=} F_{j,j+1}, \quad j = 1, 2, \ldots, n,$$

where

$$G_{j,j+1} \stackrel{\text{def}}{=} x(e_j \cdot e_{j+1}) - z_j z_{j+1} + y_{j,j+1}, \quad H_j \stackrel{\text{def}}{=} |R \times e_j|,$$

$$z_j \stackrel{\text{def}}{=} R \cdot e_j, \quad y_{j,j+1} \stackrel{\text{def}}{=} i|R|(e_j \times e_{j+1}) \cdot R.$$

Here $e_j \times e_{j+1}$ denotes the vector product and $i \in \mathbb{C}$ is the imaginary unit, $i^2 = -1$.

Now the *fundamental flexing equation* derived in [7] is

$$\prod_{j=1}^{n} F_{j,j+1} = 1. \tag{2.1}$$

In a sense, (2.1) says that the suspension stays closed up as $x$ varies.

**The roots and branching points.**

As it is shown in [7], using more analytic geometry, we can prove that

$$y_{j,j+1} = -y_{j+1,j}, \quad z_j = \frac{1}{2}(x + e_j \cdot e_j - e'_j \cdot e'_j),$$

$$G_{j,j+1} G_{j+1,j} = H_j^2 H_{j+1}^2, \quad H_j^2 = x(e_j \cdot e_j) - z_j^2 = -\frac{1}{4}(x - r'_j)(x - r_j),$$

where $r_j \stackrel{\text{def}}{=} (|e_j| + |e'_j|)^2$ and $r'_j \stackrel{\text{def}}{=} (|e_j| - |e'_j|)^2$ are the roots of $H_j^2$, and

$$y_{j,j+1}^2 = -x \det \begin{pmatrix} e_j \cdot e_j & e_j \cdot e_{j+1} & z_j \\ e_j \cdot e_{j+1} & e_{j+1} \cdot e_{j+1} & z_{j+1} \\ z_j & z_{j+1} & x \end{pmatrix}$$

$$= \frac{1}{4}(e_{j,j+1} \cdot e_{j,j+1}) x (x - b'_{j,j+1})(x - b_{j,j+1}).$$

Here $b_{j,j+1}$, $b'_{j,j+1}$ are the *branching points* of $y^2_{j,j+1}$. They can be described as the maximum and minimum value of $x = |R|^2$ for the two triangles determined by $e_j, e_{j+1}$ and $e'_j, e'_{j+1}$ which share the common edge $e_{j,j+1} = e_j - e_{j+1}$. Namely, when $x = b_{j,j+1}$ the two triangles are planar with $N, S$ in the plane on opposite sides of the line determined by $e_{j,j+1}$. Similarly, when $x = b'_{j,j+1}$, $N, S$ are on the same side of the line. In particular, $b'_{j,j+1}, b_{j,j+1}$ are real, nonnegative, and

$$0 \leqslant r'_j \leqslant b'_{j,j+1} \leqslant b_{j,j+1} \leqslant r_j.$$

We say $y_{j,j+1}$ is *equivalent* to $y_{k,k+1}$ or that $j$ is equivalent to $k$ if $b'_{j,j+1} = b'_{k,k+1}$ and $b_{j,j+1} = b_{k,k+1}$.

Now we define $\varepsilon_{j,j+1}$ as the sign of $(e_j \times e_{j+1}) \cdot R$. Note that $\varepsilon_{j,j+1}$ is also the signe of $\theta_{j,j+1}$, and if the orientation of the suspension is chosen correctly $\varepsilon_{j,j+1}$ is $+1$ or $-1$ as the suspension is locally convex or concave, respectively, at $e_{j,j+1}$. In any case, when $x$ is in the flexing interval

$$y_{j,j+1} = \frac{i}{2}\varepsilon_{j,j+1}|e_{j,j+1}|\sqrt{-x(x - b'_{j,j+1}(x - b_{j,j+1})}$$

where the positive square root is chosen.

In [7] it is shown that, studying (2.1) near the branching points $b'_{j,j+1}$ and $b_{j,j+1}$, we can split (2.1) into several equations, each corresponding to some equivalent class of $j$'s as described in the following lemma.

**Lemma.** *Let $\mathcal{P}$ be a suspension that flexes with variable $x$, and let $\mathcal{C}_0 \subset \{1, 2, \ldots, n\}$ be a subset corresponding to an equivalence class described above. Then*

$$\prod_{j \in \mathcal{C}_0}(Q_{j,j+1} + y_{j,j+1}) = \prod_{j \in \mathcal{C}_0}(Q_{j,j+1} - y_{j,j+1}) \qquad (2.2)$$

*is an identity in $x$, where $Q_{j,j+1} \stackrel{\text{def}}{=} x(e_j \cdot e_{j+1}) - z_j z_{j+1}$ (or equivalently $Q_{j,j+1} + y_{j,j+1} = G_{j,j+1}$).*

Note that if (2.2) holds for each equivalence class $\mathcal{C}_0$ then (2.1) holds. Thus, if we can construct a suspension such that (2.2) holds for each $\mathcal{C}_0$, we will have a flexor.

In order to construct a flexible suspension we have to choose edge lengths $|e_j|$, $|e'_j|$, and $|e_{j,j+1}|$ ($j = 1, 2, \ldots, n$) such that equation (2.1) (or equation (2.2) for every class $\mathcal{C}_0$) is an identity in $x$. As it was mentioned above,

$$(Q_{j,j+1} + y_{j,j+1})(Q_{j,j+1} - y_{j,j+1}) = G_{j,j+1}G_{j+1,j} = H^2_j H^2_{j+1}$$
$$= \frac{1}{16}(x - r'_j)(x - r_j)(x - r'_{j+1})(x - r_{j+1}) \qquad (2.3)$$

and the four roots in (2.3) are entirely arbitrary up to the conditions imposed on them that all the $r'_j$'s be smaller than the smallest $r_j$. Also it is easy to see that the four roots of (2.3) determine $|e_j|$ and $|e'_j|$ (but one does not know the order) by

$$\{|e_j|, |e'_j|\} = \left\{\frac{1}{2}(\sqrt{r_j} + \sqrt{r'_j}), \frac{1}{2}(\sqrt{r_j} - \sqrt{r'_j})\right\}.$$

The other parameters used to define the factors on the left of (2.3) are $b_{j,j+1}$ and $b'_{j,j+1}$, and implicitly we shall discuss their relationship to the $r_j$'s later.

We now consider a fixed $\mathcal{C}_0$, with some $b'$, $b$, and define $y = i\sqrt{-x(x-b')(x-b)}$ so that

$$y^2 = x(x-b')(x-b). \tag{2.4}$$

The roots of $G_{j,j+1}$ (i.e., of $Q_{j,j+1} + y_{j,j+1}$) are the intersections of the curve defined by (2.4) and the quadratic

$$y = \frac{2Q_{j,j+1}(x)}{\varepsilon_{j,j+1}|e_{j,j+1}|}. \tag{2.5}$$

**Angle sign edge length lemmas.**

**Lemma.** [7, Lemmas 3 and 4] *Let $\mathcal{P}$ be a flexible suspension and let $\mathcal{C}_0$ be an equivalence class of $y_{j,j+1}$'s. Then*

$$\sum_{j \in \mathcal{C}_0} \varepsilon_{j,j+1}|e_{j,j+1}| = 0 \tag{2.6}$$

*and*

$$\sum_{j \in \mathcal{C}_0} y_{j,j+1} = 0. \tag{2.7}$$

Note that in [7] it is proven that (2.7) implies that the oriented volume of a flexible suspension identically equals zero during a flex and, thus, every embedded suspension (as well as every immersed suspension bounding an immersed 3-manifold) is not flexible.

**The symmetry of the roots.**

Note $y_{j,j+1} = \frac{1}{2}\varepsilon_{j,j+1}|e_{j,j+1}|y$. Thus we may regard both sides of (2.2) as polynomials in $x$ and $y$ and a root as a pair $(x, y)$. Then (2.2) simply says that $(x_0, y_0)$ is a root of the left side say, if and only if $(x_0, -y_0)$ is also a root. This in turn says that the intersections defined by the curve (2.4) and all the curves defined by (2.5) are symmetric about the $x$-axis. Also it is not hard to see that if we have quadratics defined by (2.5) and the intersections are symmetric about the $x$-axis, then (2.2) holds.

One may be tempted into guessing that the symmetry condition implies that the quadratic factors of (2.2) cancel, but in fact this does not necessarily happen.

**The non-singular cubic.**

We now arrive to the problem of how to describe in reasonably general terms how one creates factors with the symmetry condition of above.

The non-singular cubic, of which (2.2) is an example, is an abelian variety. It turns out that it is possible to define a group operation on the curve in very natural way. Namely, we can choose any point and call it 0. We shall choose 0 to be the point at $\infty$ on the $y$-axis. Then if $Q_1, Q_2, Q_3$ are three distinct points on the intersection of a line with the curve, or two of the $Q_j$'s are equal and the line is tangent to the curve there, the group is defined by the condition $Q_1 + Q_2 + Q_3 = 0$.

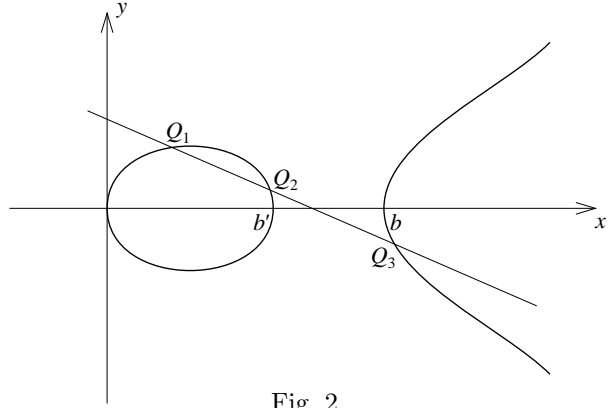

Fig. 2.

If $Q$ is on the curve, $-Q$ is the reflection of $Q$ about the $x$-axis, see Fig. 2. It is well known that this in fact defines an abelian group [17].

**The quadratic.**

Our basic problem is to describe how unsymmetric quadratics can intersect the cubic (2.4) in such a way that the intersections *are* symmetric.

Let $y = \widehat{Q}(x)$ be a quadratic curve, where $\widehat{Q}(x)$ is a quadratic function of $x$. It is easy to see that this curve intersects (2.4) at four finite points (perhaps complex points in general, but in our case they are always real). It is also easy to see that if we homoginize the equations (complete everything to projective situation) that there is in fact a double root at $\infty$ (what we called the origin before) thus giving with Bezout's theorem. Let $Q_1, Q_2, Q_3, Q_4$ be the four finite intersections of $y = \widehat{Q}(x)$ with (2.4). Then by well known results of algebraic geometry (e.g., Theorem 9.2 of [17]) we see that $Q_1 + Q_2 + Q_3 + Q_4 = 0$, and this condition is sufficient for the existence of such a $y = \widehat{Q}(x)$ to intersect (2.4) at the four given points.

**The conditions.**

We wish to write down a collection of conditions that must be satisfied if a suspension is to flex (with variable $x$). However, we need a certain amount of notation. Let
$$\widehat{Q}(x) = \frac{2Q_{j,j+1}(x)}{\varepsilon_{j,j+1}|e_{j,j+1}|}$$
be the quadratic of (2.5). We have the four roots of (2.3) which serve as the intersection of (2.4) and (2.5), and we need a way of labeling them.

We label the four points $(x,y)$ on the curve (2.4) $Q'_{j-}$, $Q_{j-}$, $Q'_{j+1+}$, $Q_{j+1+}$ corresponding respectively to the $x$-values $r'_j$, $r_j$, $r'_{j+1}$, $r_{j+1}$. Also for each point $Q$ on the curve (2.4) let $\overline{Q}$ denote its $x$ coordinate.

If the $Q_{j\pm}$'s correspond to the roots of (2.2) they must satisfy the following conditions.

(A) $\overline{Q}'_{j-} = \overline{Q}'_{j+}$, $\overline{Q}_{j-} = \overline{Q}_{j+}$, for all $j = 1, 2, \ldots, n$.

(B) $Q'_{j-} + Q_{j-} + Q'_{j+1+} + Q_{j+1+} = 0$ in the group corresponding to $j \in \mathcal{C}_0$.

(C) For every equivalence class $\mathcal{C}_0$, the collection of $Q$'s (counting multiplicities) is symmetric about the $x$-axis. (Also, the $Q'$'s are on the finite component and the $Q$'s on the infinite component.)

(B) and (C) have been discussed above. (A) is simply the condition that $r_j$ and $r'_j$ depend only on $|e_j|$ and $|e'_j|$. Also it is handy to note that if $j$ and $j+1$ are in the same equivalence class (thus defining the same curve (2.4) and the same group), then (A) is just the condition that $Q'_{j-} = \pm Q'_{j+}$, $Q_{j-} = \pm Q_{j+}$.

Using the conditions (A), (B), (C) it is possible to write down the points on a curve (2.4) that would hopefully come from a non-trivial flexor.

**The flow graph.**

The conditions above are sufficient to enable one to create many non-trivial flexors. In this subsection we describe a *flow graph* associated to a non-trivial flexor which considerably simplifies the construction of points which satisfy the conditions (A), (B), (C).

We construct a graph, $G_{\mathcal{C}_0}$, (a multi-graph in the sense of F. Harary [11]), corresponding to each equivalence class or group as follows: The vertices of $G_{\mathcal{C}_0}$ consist of the elements $j \in \mathcal{C}_0$. By property (C) there is a pairing between the roots, the $Q_{j\pm}$'s and $Q'_{j\pm}$'s. Choose one such pairing. We say $j$ is *adjacent* to $k$ if one of the $Q$'s for $j$, $Q'_{j-}$, $Q_{j-}$, $Q'_{j+1+}$, $Q_{j+1+}$, is paired with minus (in the group) one of the $Q$'s for $k$.

We furthermore wish to define a *flow* on $G_{\mathcal{C}_0}$ in the nature described in [4]. Assign a direction to the edges of $G_{\mathcal{C}_0}$ arbitrarily. If the direction of an edge is from $j$ to $k$, then the flow is $X$, if $X$ is the value (thought of as in the group) of the $j$ $Q$ that is paired with $k$. Note that condition (B) implies that the total flow into any vertex is zero.

Since each $j$ corresponds to four $Q$'s, the degree (not counting direction) at each vertex is four. Notice, also, from the nature of the graph and the fact that (2.4) has two components each $Q'$ is necessarily paired with a $Q'$ and similarly for the $Q$'s. This is because the $Q'$'s are all on the finite component of the cubic. Thus the edges can be partitioned into two equal collections, corresponding to the $Q'$'s and the $Q$'s, and each vertex is adjacent to two edges of each type. Each collection of edges is called a *two-factor* and we call them $F'$ and $F$. Thus the graph obtained is simply a graph with two disjoint two-factors, that also has a non-trivial flow.

**An example.**

Consider the flow graph

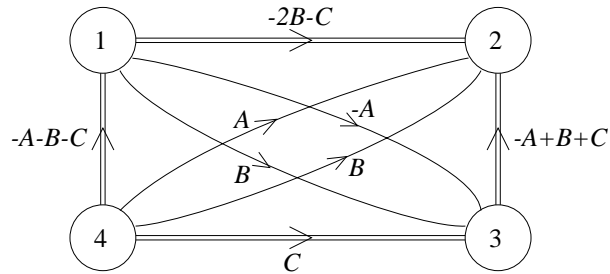

Fig. 3.

which represents the associated graph with double lines being the $F'$ factors, and single lines being the $F$ factors. The flow graph shown on Fig. 3 generates the following table

| $j$ | $Q_{j-1,-}$ | $Q_{j+}$ | $Q'_{j-1,-}$ | $Q'_{j+}$ |
|---|---|---|---|---|
| 1 | $A$ | $B$ | $C$ | $-A-B-C$ |
| 2 | $B$ | $-A$ | $A+B+C$ | $-2B-C$ |
| 3 | $-A$ | $-B$ | $2B+C$ | $A-B-C$ |
| 4 | $-B$ | $A$ | $-A+B+C$ | $-C$ |

Tab. 1.

Note that condition (A) puts additional constraints an the flow and is automatically incorporated in the above flows.

Table 1 corresponds to the first type flexible octahedra, which are constructed by taking a quadralateral in the plane that has opposite edges equal, but crosses itself, and then choosing the north and south pole in the plane of symmetry through the crossing point, see Fig. 4. These were described, e. g., in [14].

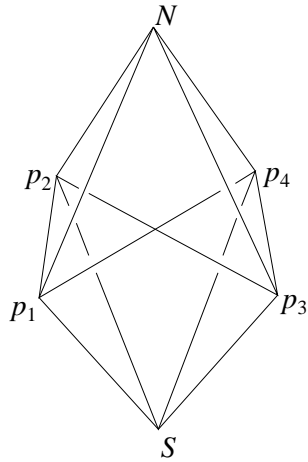

Fig. 4.

§ 3. A FLEXIBLE SUSPENSION WITH A HEXAGONAL EQUATOR

In this section we construct a flexible suspension with a hexagonal equator.

**Step 1: choosing the flow graph.**

Consider the directed multigraph shown on Fig. 5. Suppose the flows corresponding to the double lines are called $F'$ factors and the flows corresponding to the single lines are called $F$ factors.

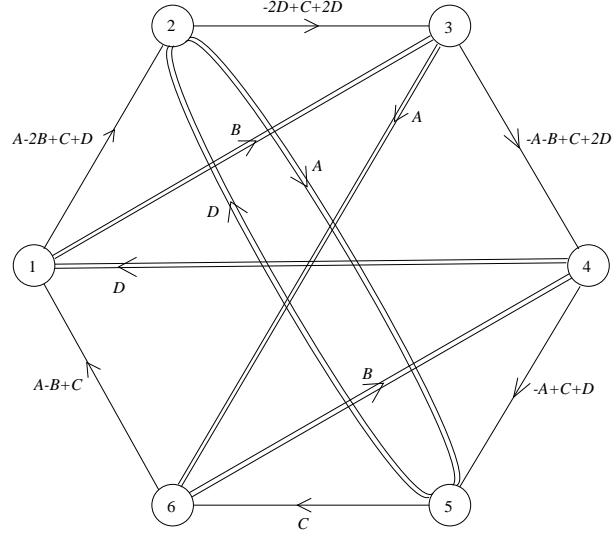

Fig. 5.

Note that this graph satisfies the conditions (A), (B), (C), mentioned in § 2.

**Step 2: calculating the table.**

Consider the following table

| $j$ | $Q_{j-1,-}$ | $Q_{j+}$ | $Q'_{j-1,-}$ | $Q'_{j+}$ |
|---|---|---|---|---|
| 1 | $C$ | $-A+B-C$ | $A$ | $-B$ |
| 2 | $A-B+C$ | $-A+2B-C-D$ | $-B$ | $D$ |
| 3 | $A-2B+C+D$ | $2B-C-2D$ | $D$ | $-A$ |
| 4 | $-2B+C+2D$ | $A+B-C-2D$ | $-A$ | $B$ |
| 5 | $-A-B+C+2D$ | $A-C-D$ | $B$ | $-D$ |
| 6 | $-A+C+D$ | $-C$ | $-D$ | $A$ |

Tab. 2.

Note that that Tab. 2 corresponds to (or is generated by) the flow graph shown on Fig. 5.

Our goal is to construct a flexible suspension with the unique equivalence class $\mathcal{C}$ and with the table shown in Tab. 2.

**Step 3: calculating the points on the cubic.**

We fix $b' = 51$, $b = 100$ and consider $y$ as an algebraic function of $x$ determined by the equation
$$y^2 = x(x-b')(x-b), \qquad (3.1)$$
see Fig. 6. By definition, we put

$$A = (2, 98); \ B = \left(\frac{4039540}{762129}, \frac{100768585960}{665338617}\right);$$
$$C = (102, -102); \ D = (30, -210). \qquad (3.2)$$

Obviously, points $A$, $B$, and $D$ lie on the bounded component of non-singular cubic (3.1) while point $C$ lies on the unbounded component, as it is shown on Fig. 6.

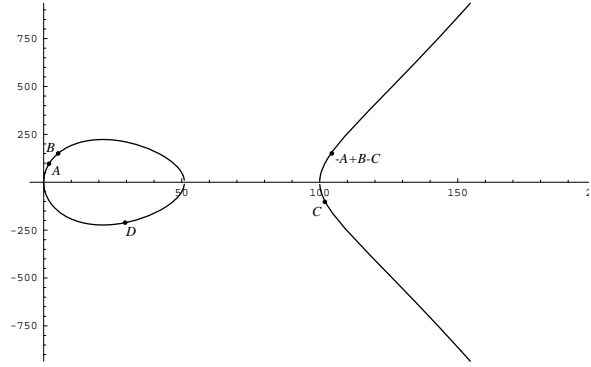

Fig. 6.

Using either the definition of the group addition on non-singular cubic (3.1) given above in §2 or an explicit formula for the group addition derived in [12] with the help of the Weirstrass function, we write the coordinates of the sum of the points $(u_1, v_1)$ and $(u_2, v_2)$ on the cubic as follows

$$(u_1, v_1) + (u_2, v_2) = \left(b' + b - u_1 - u_2 + \left(\frac{v_1 - v_2}{u_1 - u_2}\right)^2,\right.$$
$$\left. v_1 - \frac{v_1 - v_2}{u_1 - u_2}\left[b' + b - 2u_1 - u_2 + \left(\frac{v_1 - v_2}{u_1 - u_2}\right)^2\right]\right). \quad (3.3)$$

Using formula (3.3) and the system of analytic computations *Mathematica* we calculate coordinates of the other points mentioned in Table 2[1]

$$-A + B - C = \left(\frac{30931440}{292681}, \frac{28695544920}{158340421}\right);$$
$$-A + 2B - C - D = \left(\frac{56616292808335490583277 70}{39463950615542162398 09},\right.$$
$$\left.\frac{127603537646309568645683852689555598 30}{2479138777711820010358825586537 7}\right); \quad (3.4)$$
$$2B - C - 2D = \left(\frac{98365674940749318}{521862179555809}, -\frac{180534115140397956856257 54}{11921577754013306206127}\right);$$
$$A + B - C - 2D = (240, -2520);$$
$$A - C - D = \left(\frac{49130}{121}, \frac{8840510}{1331}\right).$$

Point $-A + B - C$ is also shown on Fig. 6, while the other points are too far from the origin to be shown there.

**Step 4: calculating the $r'_j$'s, $r_j$'s, and $\varepsilon_{j,j+1}$'s.**

Recall that in § 2 we label the points $Q'_{j-}$ and $Q_{j-}$ on the curve (2.4) in such a way that their $x$-values are $r'_j$ and $r_j$ respectively. Thus Table 2 and formulas (3.2) and (3.3) immediately give us the following values for the $r'_j$'s and $r_j$'s:

---

[1] According to § 2, if point $X$ has coordinates $(u, v)$ then point $-X$ has coordinates $(u, -v)$ and, thus, can be easily found. This is the reason why we write in (3.4) and show on Fig. 6 only one of the points $X$ and $-X$.

| $j$ | $r'_j$ | $r_j$ | $\varepsilon_{j,j+1}$ |
|---|---|---|---|
| 1 | $\dfrac{4039540}{762129} \approx 5.30$ | $\dfrac{30931440}{292681} \approx 105.68$ | $-1$ |
| 2 | 30 | $\dfrac{566162928083354 9058327770}{39463950615542 16239809} \approx 1434.63$ | $+1$ |
| 3 | 2 | $\dfrac{98365674940749318}{521862179555809} \approx 188.49$ | $+1$ |
| 4 | $\dfrac{4039540}{762129} \approx 5.30$ | 240 | $-1$ |
| 5 | 30 | $\dfrac{49130}{121} \approx 406.03$ | $+1$ |
| 6 | 2 | 102 | $-1$ |

Tab. 3.

Note that, using formulas (3.2) and (3.4), we can read from the appropriate line in Tab. 2 all coordinates of the four points of intersection of cubic (3.1) and the quadratic (2.5), i.e.,

$$y = \frac{2Q_{j,j+1}(x)}{\varepsilon_{j,j+1}|e_{j,j+1}|}.$$

Hence the branches of (2.5) are oriented upward if and only if the $y$-coordinate of the right-most of the four points of intersection is positive. Now, taking into account that the leading coefficient of $Q_{j,j+1}(x)$ equals $-1/4$, we decide whether $\varepsilon_{j,j+1} = +1$ or $\varepsilon_{j,j+1} = -1$.

For example,

$$Q'_1 = -B = \left(\frac{4039540}{762129}, -\frac{100768585960}{665338617}\right) \approx (5.30, -151.45),$$

$$Q'_2 = D = (30, -210),$$

$$Q_1 = A - B + C = \left(\frac{30931440}{292681}, -\frac{28695544920}{158340421}\right) \approx (105.68, -181.23),$$

$$Q_2 = -A + 2B - C - D = \left(\frac{566162928083354 9058327770}{39463950615542 16239809},\right.$$
$$\left.\frac{127603537646309568645683852689555 59830}{247913877777118200103588255865377}\right) \approx (1434.63, 51470.90).$$

The right-most point is $Q_2$ and its $y$-coordinate is positive. Hence the leading coefficient of quadratic (2.5) is positive and $\varepsilon_{1,2} = -1$.

Proceeding in the same way, we find $\varepsilon_{j,j+1}$ for $j = 2, \ldots, 6$ and put the results on Tab. 3.

**Step 5: calculating the $|e'_j|$'s, $|e_j|$'s, and $|e_{ij}|$'s.**

Recall from § 2 that the lengths of the edges of the suspension adjacent to the north and south poles, $|e_j|$ and $|e'_j|$ respectively, are such that

$$\{|e_j|, |e'_j|\} = \left\{\frac{1}{2}(\sqrt{r_j} + \sqrt{r'_j}), \frac{1}{2}(\sqrt{r_j} - \sqrt{r'_j})\right\}. \tag{3.5}$$

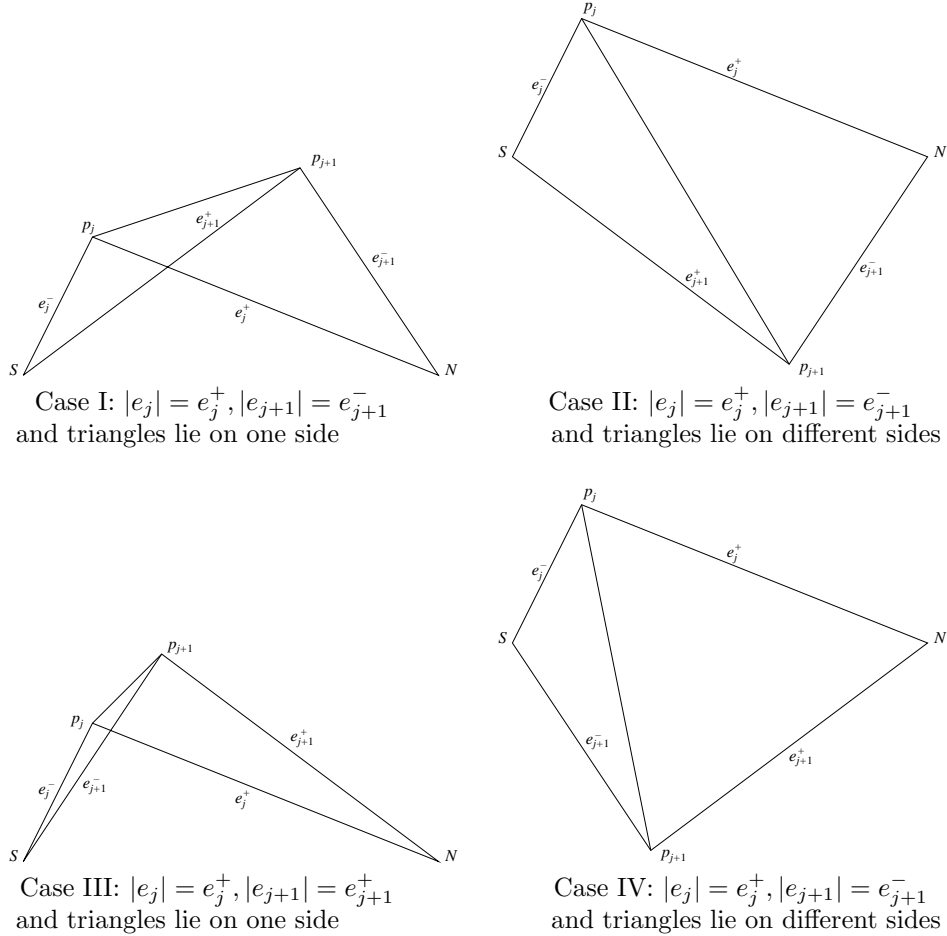

Fig. 7

Note that this relation does not allow us to determine $|e_j|$ and $|e'_j|$ precisely because we do not know the order.

Let us call the union of triangles $\langle N, p_j, p_{j+1}\rangle$ and $\langle S, p_j, p_{j+1}\rangle$ the $j$-sector of the suspension. Obviously, every $j$-sector has two flat positions. In one flat position the both triangles lie on the same side of the line determined by the equator edge $\langle p_j, p_{j+1}\rangle$ while in the other flat position the triangles lie on the different sides of that line. Using the notation $e_j^- = \min\{|e'_j|, |e_j|\}$ and $e_j^+ = \max\{|e'_j|, |e_j|\}$ we draw four possible cases of mutual location of triangles $\langle N, p_j, p_{j+1}\rangle$ and $\langle S, p_j, p_{j+1}\rangle$ on Fig. 7. The four cases shown on Fig. 7 correspond to the situation when $|e_j| = e_j^+$. The other four cases, corresponding to the situation when $|e_j| = e_j^-$, can be drown similarly. We do not draw those four cases because they can be obtained from the cases shown on Fig. 7 by interchanging the north and south poles, i. e., by replacing $S$ by $N$ and $N$ by $S$ on Fig. 7.

Consider Case I shown on Fig. 7. Observe that the same sector has yet another flat position which may be obtained by rotating triangle $\langle N, p_j, p_{j+1}\rangle$ around its side $\langle p_j, p_{j+1}\rangle$ in 3-space to the angle $\pi$. The result is shown on Fig. 8. Now triangles $\langle N, p_j, p_{j+1}\rangle$ and $\langle S, p_j, p_{j+1}\rangle$ lie on one side of the line through vertices $p_j$ and $p_{j+1}$. We say that the configuration shown on Fig. 8 is obtained from Case I on Fig. 7 by a *flip*. Similarly, we may apply a flip to each of the Cases II–IV shown on Fig. 7.

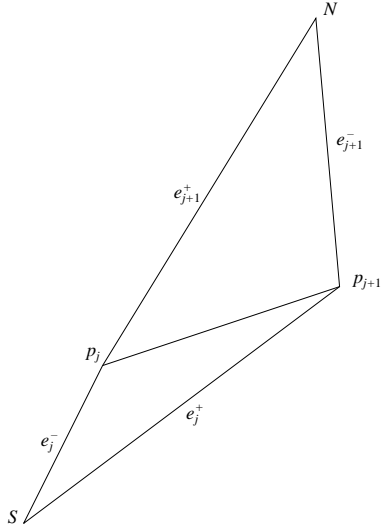

Fig. 8.

For every sector we know from §2 that, in one flat position, the distance between the north and south poles, $N$ and $S$, equals $b$, while in the other flat position it is equal to $b'$. At the very beginning of §3 we have fixed $b' = 51$ and $b = 100$.

For each $j = 1, \ldots, 6$ we use Table 3 to compute

$$e_j^- = \min\left\{\frac{1}{2}(\sqrt{r_j} + \sqrt{r'_j}), \frac{1}{2}(\sqrt{r_j} - \sqrt{r'_j})\right\}$$

and

$$e_j^+ = \max\left\{\frac{1}{2}(\sqrt{r_j} + \sqrt{r'_j}), \frac{1}{2}(\sqrt{r_j} - \sqrt{r'_j})\right\}.$$

Then, for each of the Cases I–IV, we put the poles, $S$ and $N$, at the points $(0,0)$ and $(10,0)^2$ respectively and calculate the coordinates of the points $p_j$ and $p_{j+1}$. The distance $|p_j - p_{j+1}|$ is a candidate for the length $e_{j,j+1}$ of the equatorial edge. At last, we calculate the coordinates of the north pole under the flip transformation (i.e., under the reflection of triangle $\langle N, p_j, p_{j+1}\rangle$ with respect to the line through the points $p_j$ and $p_{j+1}$) and if this distance equals exactly $\sqrt{51}$ (i.e., this distance squared equals $b' = 51$ that is the square of the minimal distance between the poles during the flex) then we say that the choice for the $e_j$, $e'_j$, $e_{j+1}$, $e'_{j+1}$, and $e_{j,j+1}$ is correct. Note that, when we calculate $e_j$, $e'_j$, $e_{j+1}$, $e'_{j+1}$, and $e_{j,j+1}$ for the first sector we can obviously interchange the north and south poles, but as soon as the poles, $N$ and $S$, are fixed for the first sector, the above described procedure determine the edge lengths in a unique way. We accumulate the results of such calculations in Tables 4 and 5.

---

[2]Recall from §2 that $b = 100$ is the square of the maximal distance between the poles during the flex.

| $j$ | $|e_j|$ | $|e'_j|$ |
|---|---|---|
| 1 | $\frac{1}{2}\left(\frac{1436\sqrt{15}}{541} + \frac{218\sqrt{85}}{873}\right) \approx 6.29$ | $\frac{1}{2}\left(\frac{1436\sqrt{15}}{541} - \frac{218\sqrt{85}}{873}\right) \approx 3.99$ |
| 2 | $\frac{1}{2}\left(\frac{182493018091\sqrt{170}}{62820339553} - \sqrt{30}\right) \approx 16.20$ | $\frac{1}{2}\left(\frac{182493018091\sqrt{170}}{62820339553} + \sqrt{30}\right) \approx 21.68$ |
| 3 | $\frac{1}{2}\left(\frac{31054297\sqrt{102}}{22844303} - \sqrt{2}\right) \approx 6.16$ | $\frac{1}{2}\left(\frac{31054297\sqrt{102}}{22844303} + \sqrt{2}\right) \approx 7.57$ |
| 4 | $\frac{1}{2}\left(4\sqrt{15} - \frac{218\sqrt{85}}{873}\right) \approx 6.59$ | $\frac{1}{2}\left(4\sqrt{15} + \frac{218\sqrt{85}}{873}\right) \approx 8.90$ |
| 5 | $\frac{1}{2}\left(\frac{17\sqrt{170}}{11} - \sqrt{30}\right) \approx 7.34$ | $\frac{1}{2}\left(\frac{17\sqrt{170}}{11} + \sqrt{30}\right) \approx 12.81$ |
| 6 | $\frac{1}{2}(\sqrt{102} + \sqrt{2}) \approx 5.76$ | $\frac{1}{2}(\sqrt{102} - \sqrt{2}) \approx 4.34$ |

Tab. 4.

| $j$ | $|e_{j,j+1}|$ |
|---|---|
| 1 | $\frac{541419683182996345}{296699606628505029} \approx 18.25$ |
| 2 | $\frac{31635727886833754300}{1435086871311616559} \approx 22.04$ |
| 3 | $\frac{27288800741}{19943076519} \approx 1.37$ |
| 4 | $\frac{130585}{9603} \approx 13.60$ |
| 5 | $\frac{100}{11} \approx 9.09$ |
| 6 | $\frac{310327}{472293} \approx 0.66$ |

Tab. 5.

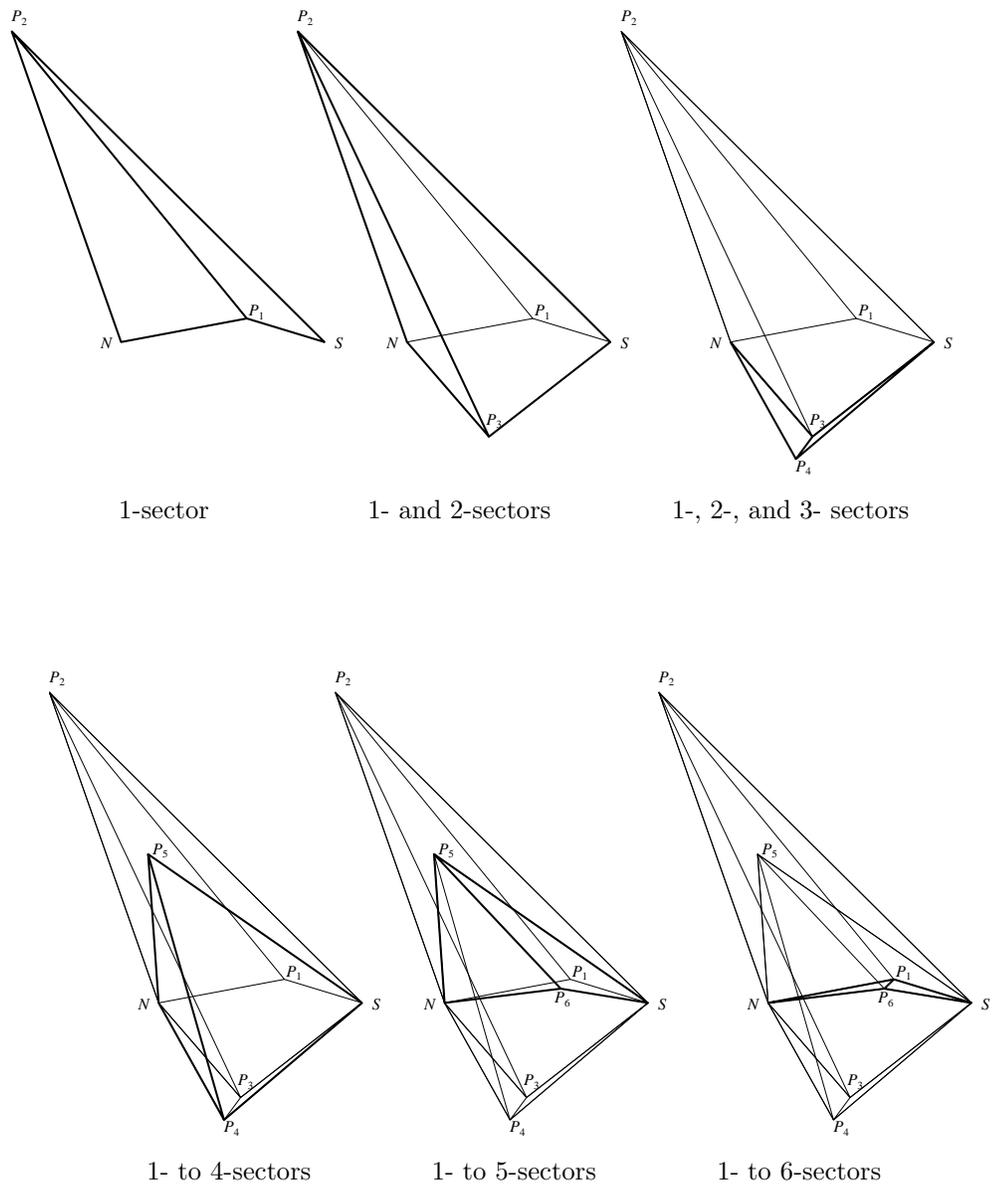

Fig. 9.

Let us mention that when we calculate $e_{61}$ we additionally take into account that, at that moment, we already know the order how the edges $e_1$, $e'_1$, $e_6$ and $e'_6$ are attached to the north and south poles, $N$ and $S$ respectively, and, of course, their lengths which are shown in Tab. 4.

On Fig. 9 we show how the next sector is glued to the previous one in a flat position corresponding to $x = 100$. To this end we draw the last glued sector with thicken lines in contrast with all the preceding sectors. It is quiet expectable that, in each sector, the triangles $\langle N, p_j, p_{j+1}\rangle$ and $\langle S, p_j, p_{j+1}\rangle$ are located on different sides of the line through the points $p_j$ and $p_{j+1}$.

Fix an orientation of the whole suspension constructed. By $\varphi_j$ denote the dihedral angle of the suspension at the edge $\langle N, p_j\rangle$ (for immersed suspensions the notion of dihedral angle is clear; for non immersed suspensions we discuss this notion in §4 below). Similarly, denote by $\varphi'_j$ and $\varphi_{j,j+1}$ the dihedral angles at the edges $\langle S, p'_j\rangle$ and $\langle p_j, p_{j+1}\rangle$. Note that these angles are functions in $x$, the squared distance between the north and south poles, $N$ and $S$. Obviously, Fig. 9 provides us with a possibility to find those angles modulo $2\pi$ for $x = 100$. The results are presented in Table 6.

| $j$ | $\varphi_j(100)$ | $\varphi'_j(100)$ | $\varphi_{j,j+1}(100)$ | $\varphi_j(51)$ | $\varphi'_j(51)$ | $\varphi_{j,j+1}(51)$ |
|---|---|---|---|---|---|---|
| 1 | $\pi$ | $\pi$ | $\pi$ | 0 | 0 | 0 |
| 2 | 0 | 0 | $\pi$ | $\pi$ | $\pi$ | 0 |
| 3 | $\pi$ | $\pi$ | $\pi$ | 0 | 0 | 0 |
| 4 | 0 | 0 | $\pi$ | $\pi$ | $\pi$ | 0 |
| 5 | 0 | 0 | $\pi$ | $\pi$ | $\pi$ | 0 |
| 6 | 0 | 0 | $\pi$ | $\pi$ | $\pi$ | 0 |

Tab. 6.

On Fig. 10 we show how the next sector is glued to the previous one in a flat position corresponding to $x = 51$. As before we draw the last glued sector with thicken lines in contrast with all the preceding sectors. It is quiet expectable that, in each sector, the triangles $\langle N, p_j, p_{j+1}\rangle$ and $\langle S, p_j, p_{j+1}\rangle$ are located on the same side of the line through the points $p_j$ and $p_{j+1}$.

Obviously, Fig. 10 provides us with a possibility to find those angles modulo $2\pi$ for $x = 51$. The results are presented in Table 6.

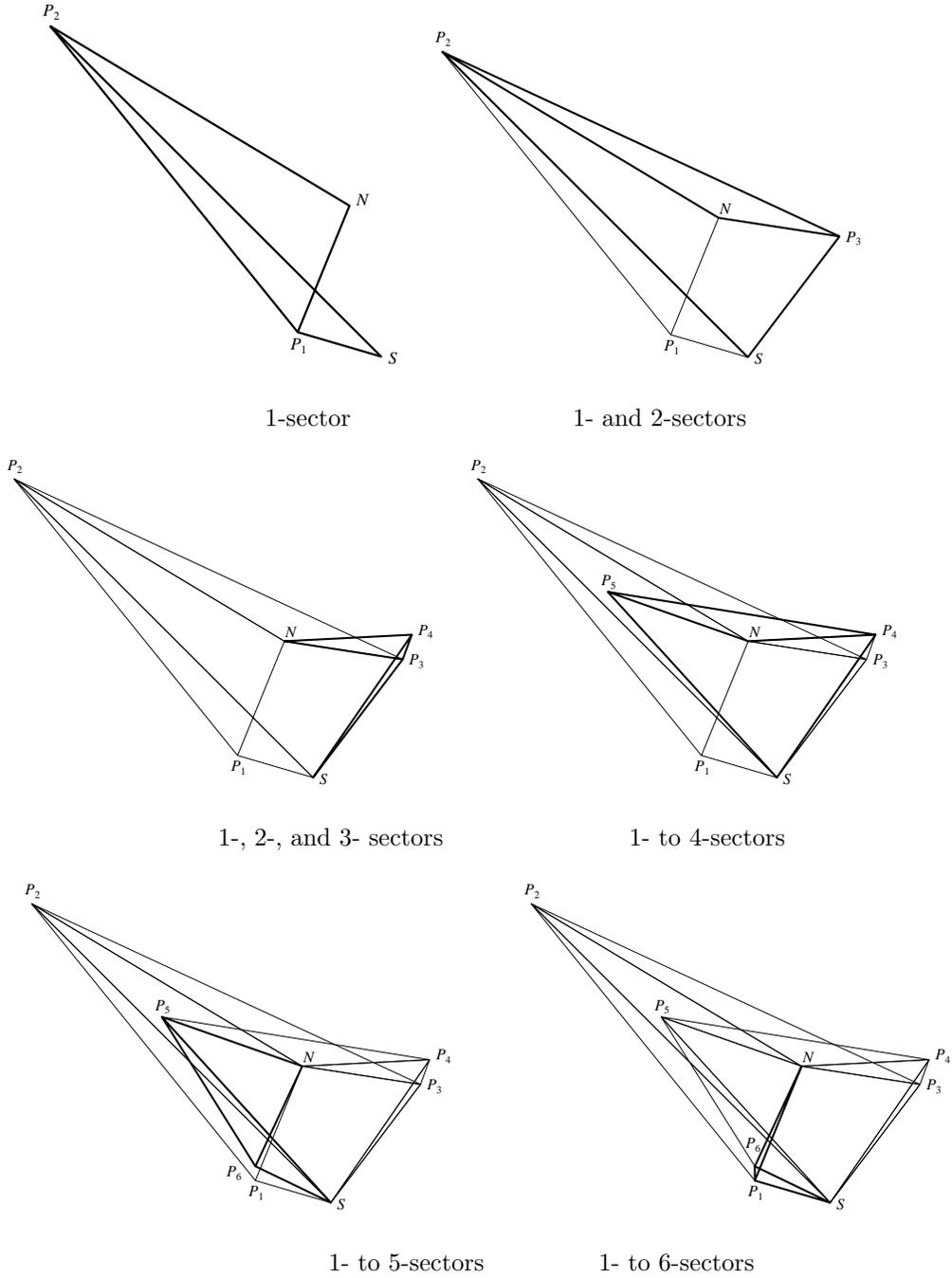

Fig. 10.

**Comparison with previously known theorems.**

It follows from a Mikhalev's theorem [13] that if a suspension with an $n$-gone equator $p_1, \ldots, p_n$ is bent in such a way that the length of the 'short' diagonal $\langle p_{j-1}, p_{j+1} \rangle$ is non-constant for some $j = 1, \ldots, n$ then there exists $k = 1, \ldots, n$, $k \neq j \pm 1$, such that

$$|e_j| + (-1)^{\sigma_1}|e'_j| + (-1)^{\sigma_2}|e_k| + (-1)^{\sigma_3}|e'_k| = 0 \tag{3.6}$$

with some integers $\sigma_1$, $\sigma_2$, $\sigma_3$.

Note that (3.6) is satisfied for the above constructed suspension $\mathcal{S}$ with a hexagonal equator, namely,

$$|e_1| - |e'_1| + |e_4| - |e'_4| = 0,$$
$$|e_2| - |e'_2| - |e_5| + |e'_5| = 0,$$
$$|e_3| - |e'_3| + |e_6| - |e'_6| = 0.$$

It follows from another Mikhalev's theorem [13] that, for every flexible suspension with an $n$-gone equator $p_1, \ldots, p_n$

$$\sum_{j=1}^{n} (-1)^{\sigma_j}|e_{j,j+1}| = 0 \tag{3.7}$$

with some integers $\sigma_j$.

Note that (3.7) obviously follows from (2.6) while the latter was proven in [7] more than 25 years prior to Mikhalev's paper [13]. Moreover, (2.6) provides us with additional geometric information that (3.7) holds true with $(-1)^{\sigma_j} = \varepsilon_{j,j+1}$, where, as it was specified in § 2, $\varepsilon_{j,j+1}$ equals the sign of $(e_j \times e_{j+1}) \cdot R$.

Using data from Tables 3 and 5, we find by direct calculations that (3.7), or (2.6), is satisfied for the above constructed suspension $\mathcal{S}$ with a hexagonal equator, namely,

$$-|e_{12}| + |e_{23}| + |e_{34}| - |e_{45}| + |e_{56}| - |e_{61}| = 0 \tag{3.8}$$

or

$$-\frac{541419683182996345}{29669606628505029} + \frac{31635727886833754300}{1435086871311616559}$$
$$+ \frac{27288800741}{19943076519} - \frac{130585}{9603} + \frac{100}{11} - \frac{310327}{472293} = 0.$$

**Discussion.**

Obviously, every flexible octahedron (known also as a Bricard's octanedron) can be considered as a flexible suspension with a quadrilateral equator and gives rise to trivial flexible suspensions with, say, pentagonal or hexagonal equators which can be constructed as follows: fix an interior point on an equator edge and join it with the north and south poles with new edges (i.e., subdivide some faces of the octahedron).

There is another obvious way to construct a trivial flexible suspension with pentagonal or hexagonal equator: start with an arbitrary suspension with a pentagonal or hexagonal equator; remove the star of the south pole and treat the star of the north pole as a twice-covered polyhedral surface.

The above mentioned trivial flexible suspensions, definitely, cannot help us to construct a counterexample to the Strong Bellows Conjecture. We were not able to construct a non-trivial flexible suspension with pentagonal equator, but we can summarize the results of this section in the following theorem.

**Theorem.** *The suspension $\mathcal{S}$, constructed above, provides us with a non-trivial example of a flexible suspension with a hexagonal equator.*

### §4. An attack on the Strong Bellows Conjecture

In this section we study properties of the flexible suspension with a hexagonal equator $\mathcal{S}$ constructed in § 3 which are related to the Strong Bellows Conjecture.

**The Strong Bellows Conjecture.**

In a comment added to the Russian translation of [7] R. Connelly conjectured that *if an immersed polyhedron $\mathcal{P}$ in Euclidean 3-space is obtained from another immersed polyhedron $\mathcal{Q}$ by a continuous flex then $\mathcal{P}$ and $\mathcal{Q}$ are scissors congruent*, i.e., $\mathcal{P}$ can be divided into a finite set of tetrahedra each of which can be moved independently one from another in 3-space in such a way that the resulting set constitutes a partition of $\mathcal{Q}$. The conjecture still remains open and is know as the *Strong Bellows Conjecture*.

**Dehn invariants and the Extended Strong Bellows Conjecture.**

Let $f : \mathbb{R} \to \mathbb{R}$ be a $\mathbb{Q}$-linear function such that $f(\pi) = 0$, i.e., let $f(px + qy) = pf(x) + qf(y)$ for all $p, q \in \mathbb{Q}$, $x, y \in \mathbb{R}$ and $f(\pi) = 0$. The sum

$$D_f(\mathcal{P}) = \sum f(\varphi_j) \ell_j$$

is called the *Dehn invariant* of an immersed pohyhedron $\mathcal{P}$ in Euclidean 3-space. Here $\varphi_j$ is the (internal) dihedral angle at the $j$'s edge, $\ell_j$ is the length of the $j$'s edge, and the sum is taken over all the edges of $\mathcal{P}$.

It is well-known that *two immersed polyhedra in Euclidean 3-space are scissors congruent if and only if they have the same volume and every Dehn invariant takes the same value for those polyhedra*, see [6] or [10].

It seems natural to have this theorem in mind when approaching the Strong Bellows Conjecture but the first problem here is that we should extend the notions of volume, internal dihedral angle, and Dehn invariant onto an arbitrary oriented (not necessarily immersed) polyhedron.

The extension of the notion of volume we need is the standard notion of the oriented volume [7]. In 1996, I.Kh. Sabitov [15] has proved that *every oriented flexible polyhedron in Euclidean 3-space preserves its oriented volume during a flex* and, thus, gave an affirmative answer to the Bellows Conjecture. An improved presentation is given in [16]; another proof is published in [9].

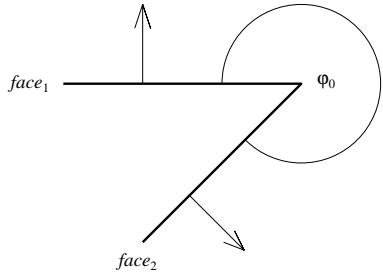

Fig. 11.

We define the *dihedral angle* at an edge $\ell$ of an oriented (not necessarily immersed) polyhedron in Euclidean 3-space as a multi-valued function $\varphi = \varphi_0 + 2\pi k$, where $k$ is an integer and one of the values, $\varphi_0$, is defined as the 'usual' dihedral angle between the two faces $f_1$ and $f_2$ adjacent to $\ell$ measured from the side determined by the orientation, see Fig. 11. In other words we can say that $\varphi_0$ is a number between 0 and $2\pi$ which is obtained as the product of $2\pi$ and the proportion of a sufficiently small ball centered at a relative interior point of $\ell$ which is contained in the intersection of the two half-spaces determined by oriented faces $f_1$ and $f_2$.

Similarly to complex analysis, we say that values of dihedral angle $\varphi$ which correspond to different values of $k$ represent different *branches* of that multi-valued function. For a flexible polyhedron we fix a particular value (or a branch) of its every dihedral angle (i.e., fix all $k$'s) in a single position and assume that those particular values change continuously in the course of the flex. This means that the value of a dihedral angle may drift from one branch to another during the flex.

The sum
$$D_f(\mathcal{P}) = \sum f(\varphi_j)\ell_j$$
is called the *Dehn invariant* of an oriented (not necessarily immersed) pohyhedron $\mathcal{P}$ in Euclidean 3-space. Here $f : \mathbb{R} \to \mathbb{R}$ is a $\mathbb{Q}$-linear function such that $f(\pi) = 0$, $\varphi_j$ is the dihedral angle at the $j$'s edge, $\ell_j$ is the length of the $j$'s edge, and the sum is taken over all the edges of $\mathcal{P}$. Note that $D_f$ does not depend on a choice of a branch of $\varphi_j$.

Now we can formulate the *Extended Strong Bellows Conjecture*: Every Dehn invariant of an oriented (not necessarily immersed) flexible polyhedron in Euclidaen 3-space remains constant during the flex.

Since we know that (oriented) volume is preserved during the flex, the Extended Strong Bellows Conjecture is equivalent to the Strong Bellows Conjecture for immersed polyhedra.

**An attack on the Extended Strong Bellows Conjecture.**

Let us treat $\mathbb{R}$ as an (infinite-dimensional) vector space over $\mathbb{Q}$. It is known that, if Zermelo's axiom is true, there exist a basis $\mathfrak{H}$ in that vector space [6] (this means that every real number $x$ is expressible uniquely in the form of a finite linear combination of elements of $\mathfrak{H}$ with rational coefficients $x = \alpha_1 e_1 + \cdots + \alpha_{n(x)} e_{n(x)}$, $\alpha_j \in \mathbb{Q}$, $e_j \in \mathfrak{H}$, $j = 1, \ldots, n(x)$, $n(x) \in \mathbb{N}$). Such a basis is known as a *Hamel base*. It is known that, without loss of generality, we may assume that $\pi$ (i.e., the area of a unit disk in Euclidean 2-space) is an element of $\mathfrak{H}$. If we assume this, we can at once write down all $\mathbb{Q}$-linear functions $f : \mathbb{R} \to \mathbb{R}$ such that $f(\pi) = 0$; we put $f(\pi) = 0$, give $f(e)$ arbitrary values for $e \in \mathfrak{H}$, $e \neq \pi$, and define $f(x)$ generally by $f(x) = \alpha_1 f(e_1) + \cdots + \alpha_n f(e_n)$ for $x = \alpha_1 e_1 + \cdots + \alpha_{n(x)} e_{n(x)}$, $\alpha_j \in \mathbb{Q}$, $e_j \in \mathfrak{H}$, $j = 1, \ldots, n(x)$. Obviously, we can represent an arbitrary $\mathbb{Q}$-linear function $f : \mathbb{R} \to \mathbb{R}$ such that $f(\pi) = 0$ as $f(x) = \sum p_e f_e(x)$, where the sum is taken over all $e \in \mathfrak{H}$, except $\pi$; $p_e$ are arbitrary real numbers; and $f_e(x) = \alpha_k$ provided $x = \alpha_1 e_1 + \cdots + \alpha_{n(x)} e_{n(x)}$ and $e = e_k$. The latter representation make it clear that the following two statements are equivalent:

(1) *Dehn invariant $D_f$ remains constant during a flex for every $\mathbb{Q}$-linear function $f : \mathbb{R} \to \mathbb{R}$ such that $f(\pi) = 0$;*

(2) *Dehn invariant $D_{f_e}$ remains constant during a flex for every dual element $f_e$, corresponding to $e \in \mathfrak{H}$, $e \neq \pi$.*

Note also that $f_e(x)$ is a rational number for every $x \in \mathbb{R}$ and every $e \in \mathfrak{H}$.

Let $\mathcal{S}$ be a flexible suspension constructed in § 3. Fix some $e \in \mathfrak{H}$, $e \neq \pi$. Substituting the values of the edge lengths, $|e_j|$, $|e'_j|$, and $|e_{j,j+1}|$, from Tables 4 and 5 to the expression of Dehn invariant

$$D_{f_e} = \sum_{j=1}^{6} \left[ f_e(\varphi_j(x))|e_j| + f_e(\varphi'_j(x))|e'_j| + f_e(\varphi_{j,j+1}(x))|e_{j,j+1}| \right]$$

yields

$$\begin{aligned} D_{f_e} =& f_e(\alpha_1(x)) + f_e(\alpha_2(x))\sqrt{2} + f_e(\alpha_3(x))\sqrt{15} \\ &+ f_e(\alpha_4(x))\sqrt{30} + f_e(\alpha_5(x))\sqrt{85} + f_e(\alpha_6(x))\sqrt{102} + f_e(\alpha_7(x))\sqrt{170}, \end{aligned}$$

where

$$\begin{aligned} \alpha_1(x) =& \frac{5414196831829963445}{29669606628505029}\varphi_{12}(x) + \frac{31635727886833754300}{14350868713111616559}\varphi_{23}(x) \\ &+ \frac{27288800741}{19943076519}\varphi_{34}(x) + \frac{130585}{9603}\varphi_{45}(x) \\ &+ \frac{100}{11}\varphi_{56}(x) + \frac{310327}{472293}\varphi_{61}(x), \end{aligned} \quad (4.1)$$

$$\alpha_2(x) = -\frac{1}{2}(\varphi_3(x) - \varphi'_3(x)) + \frac{1}{2}(\varphi_6(x) - \varphi'_6(x)), \quad (4.2)$$

$$\alpha_3(x) = \frac{718}{541}(\varphi_1(x) + \varphi'_1(x)) + 2(\varphi_4(x) + \varphi'_4(x)), \quad (4.3)$$

$$\alpha_4(x) = -\frac{1}{2}(\varphi_2(x) - \varphi'_2(x)) - \frac{1}{2}(\varphi_5(x) - \varphi'_5(x)), \quad (4.4)$$

$$\alpha_5(x) = \frac{109}{873}(\varphi_1(x) - \varphi'_1(x)) - \frac{109}{873}(\varphi_4(x) - \varphi'_4(x)), \quad (4.5)$$

$$\alpha_6(x) = \frac{31054297}{45688606}(\varphi_3(x) + \varphi'_3(x)) + \frac{1}{2}(\varphi_6(x) + \varphi'_6(x)), \quad (4.6)$$

$$\alpha_7(x) = \frac{182493018091}{125640679106}(\varphi_2(x) + \varphi'_2(x)) + \frac{17}{22}(\varphi_5(x) + \varphi'_5(x)). \quad (4.7)$$

Since the numbers 2, $15 = 3 \cdot 5$, $30 = 2 \cdot 3 \cdot 5$, $85 = 5 \cdot 17$, $102 = 2 \cdot 3 \cdot 17$, and $170 = 2 \cdot 5 \cdot 17$ are square free it follows that the numbers 1, $\sqrt{2}$, $\sqrt{15}$, $\sqrt{30}$, $\sqrt{85}$, $\sqrt{102}$, and $\sqrt{170}$ are linearly independent over rationals. Taking into consideration that $f_e(\alpha_j(x))$ is rational for all $x$ and $j = 1, \ldots, 7$, we conclude that $D_{f_e}$ is constant in $x$ if and only if $\alpha_j(x)$ is constant in $x$ for every $j = 1, \ldots, 7$.

The rest of this article is devoted to the study of expressions $\alpha_j(x)$.

**Dihedral angles adjacent to $p_2$, $p_4$, $p_5$, and $p_6$.**

By $\mathcal{Q}_j$, $j = 1, \ldots, 6$, denote the intersection of suspension $\mathcal{S}$ with a sphere centered at $p_j$ of a so small radius that it contains no vertices of $\mathcal{S}$ other than $p_j$. Note that, since $\mathcal{S}$ admits two flat positions (which occur when the distance between the poles is 10 or $\sqrt{51}$), spherical quadrangle $\mathcal{Q}_j$ admits two 'line' positions, i.e., positions when all its vertices are contained in a great circle. From Fig. 9 and Fig. 10 we conclude that, for $j = 2, 4, 5$ and 6, non of the 'line' positions of $\mathcal{Q}_j$ coincides with a whole great circle. If we denote the side lengths of $\mathcal{Q}_j$ by $\omega_{j1}, \omega_{j2}, \omega_{j3}$, and $\omega_{j4}$ in cyclic order then $\mathcal{Q}_j$ can be drawn in one 'line' position as it is shown on Fig. 12 and in another 'line' position as it is shown on Fig. 13. It follows from Fig. 12 that $\omega_{j1} + \omega_{j2} = \omega_{j3} + \omega_{j4}$. Similarly, it follows from Fig. 13 that $\omega_{j4} + \omega_{j1} = \omega_{j2} + \omega_{j3}$. Solving these two linear equations we get $\omega_{j2} = \omega_{j4}$ and $\omega_{j1} = \omega_{j3}$.

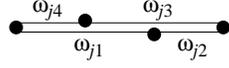
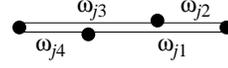

Fig. 12.          Fig. 13.

In other words, the latter means that, for $j = 2, 4, 5$, and 6, the opposite sides of $\mathcal{Q}_j$ are parewise equal to each other. Under these conditions, $\mathcal{Q}_j$ may be either a convex centrally symmetric quadrangle, see Fig. 14, or a non-convex quadrangle (with a point of self-intersection) symmetric with respect to a 'line', i.e., to a great circle, see Fig. 15. On Fig. 14 and Fig. 15 we denote the intersection of edge $\langle p_j, X \rangle$ with the small sphere centered at $p_j$ by $\widetilde{X}$.

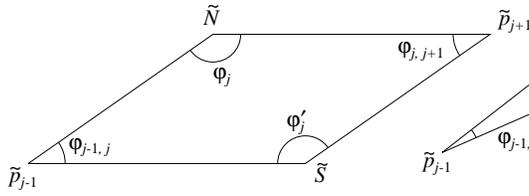
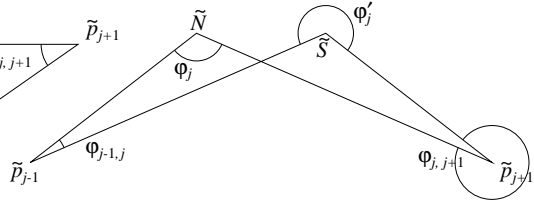

Fig. 14.          Fig. 15.

Recall from §2 that, by definition, $\varepsilon_{j,j+1}$ equals the sign of $(e_j \times e_{j+1}) \cdot R$. Equivalently, $\varepsilon_{j,j+1}$ equals the sign of $\theta_{j,j+1}$. This means that if $\varepsilon_{j-1,j}\varepsilon_{j,j+1} < 0$ then tetrahedra $\langle N, p_{j-1}, p_j, S \rangle$ and $\langle N, p_j, p_{j+1}, S \rangle$ lie on the same side of the plane which passes through the points $N$, $p_j$, and $S$. In this case $\mathcal{Q}_j$ is a non-convex spherical quadrangle whose opposite sides are pare-wise equal to each other as show on Fig. 15 and the sum of opposite angles of $\mathcal{Q}_j$ equals $2\pi$. In terms of multi-valued functions, we may write

$$\varphi'_j(x) = -\varphi_j(x) \quad \text{and} \quad \varphi_{j-1,j}(x) = -\varphi_{j,j+1}(x) \quad \text{for all } x. \tag{4.8}$$

From Tab. 3 we read $\varepsilon_{12} = \varepsilon_{45} = \varepsilon_{61} = -1$ and $\varepsilon_{23} = \varepsilon_{34} = \varepsilon_{56} = +1$. Thus, $\varepsilon_{12}\varepsilon_{23} = \varepsilon_{34}\varepsilon_{45} = \varepsilon_{45}\varepsilon_{56} = \varepsilon_{56}\varepsilon_{61} = -1 < 0$ and $\mathcal{Q}_j$ is a non-convex quadrilateral for $j = 2, 4, 5$ and 6, as described above and (4.8) hold true for the same $j$'s.

Let us summaries the relations obtained as follows

$$\begin{aligned} \varphi_{23}(x) &= -\varphi_{12}(x), & \varphi'_2(x) &= -\varphi_2(x), \\ \varphi_{45}(x) &= -\varphi_{34}(x), & \quad \text{and} \quad \varphi'_4(x) &= -\varphi_4(x), \\ \varphi_{56}(x) &= -\varphi_{45}(x), & \varphi'_5(x) &= -\varphi_5(x), \\ \varphi_{61}(x) &= -\varphi_{56}(x), & \varphi'_6(x) &= -\varphi_6(x). \end{aligned} \tag{4.9}$$

**Dihedral angles adjacent to $p_1$ and $p_3$.**

Unfortunately, we cannot apply the arguments from the previous subsection to $\mathcal{Q}_1$ and $\mathcal{Q}_3$, because in the flat position shown on Fig. 9 each of these spherical quadrangles coincides with a great circle. This is the reason why we use other arguments in this subsection.

Using the Euclidean Cosine Law and the exact values of side lengths of triangles $\langle N, p_1, p_2\rangle$, $\langle S, p_1, p_2\rangle$, $\langle S, p_1, p_6\rangle$, and $\langle N, p_1, p_6\rangle$ given in Tables 4 and 5, we find

$$\cos\angle Np_1p_2 = \frac{7(-3200524333319476 + 9569278607860305\sqrt{51})}{319877868986(626814\sqrt{15} + 58969\sqrt{85})};$$

$$\cos\angle Sp_1p_2 = -\frac{7(3200524333319476 + 9569278607860305\sqrt{51})}{319877868986(626814\sqrt{15} - 58969\sqrt{85})};$$

$$\cos\angle Sp_1p_6 = \frac{1565240 - 472293\sqrt{51}}{1253628\sqrt{15} - 117938\sqrt{85}};$$

$$\cos\angle Np_1p_6 = \frac{1565240 + 472293\sqrt{51}}{1253628\sqrt{15} + 117938\sqrt{85}}.$$

Now direct calculations show that $\arccos\angle Np_1p_2 + \arccos\angle Sp_1p_6 = \arccos\angle Sp_1p_2 + \arccos\angle Np_1p_6 = 0$. Hence

$$\angle Np_1p_2 + \angle Sp_1p_6 = \angle Sp_1p_2 + \angle Np_1p_6 = \pi. \tag{4.10}$$

Consider spherical quadrangle $\mathcal{Q}_1 = \langle \widetilde{N}, \widetilde{p_2}, \widetilde{S}, \widetilde{p_6}\rangle$ as being composed of two spherical triangles $\langle \widetilde{N}, \widetilde{S}, \widetilde{p_6}\rangle$ and $\langle \widetilde{N}, \widetilde{S}, \widetilde{p_2}\rangle$. Using the Spherical Cosine Law [5, Theorem 18.6.8] and (4.10) we get

$$\cos\angle\widetilde{p_2} = \frac{\cos\widetilde{N}\widetilde{S} - \cos\widetilde{S}\widetilde{p_2}\cos\widetilde{N}\widetilde{p_2}}{\sin\widetilde{S}\widetilde{p_2}\sin\widetilde{N}\widetilde{p_2}} = \frac{\cos\widetilde{N}\widetilde{S} - \cos\widetilde{S}\widetilde{p_6}\cos\widetilde{N}\widetilde{p_6}}{\sin\widetilde{S}\widetilde{p_6}\sin\widetilde{N}\widetilde{p_6}} = \cos\angle\widetilde{p_6},$$

where $\widetilde{N}\widetilde{S}$ stands for the spherical distance between points $\widetilde{N}$ and $\widetilde{S}$ and $\angle\widetilde{p_2}$ stands for the angle of $\mathcal{Q}_1$ at vertex $\widetilde{p_2}$.

Note that some branch of multi-valued function $\varphi_{12}(x)$ equals $\angle\widetilde{p_2}$ or $2\pi - \angle\widetilde{p_2}$ for $\mathcal{Q}_1$ constructed for the same value of $x$. Similarly, some branch of $\varphi_{61}(x)$ equals $\angle\widetilde{p_6}$ or $2\pi - \angle\widetilde{p_6}$. Taking into account that $\varepsilon_{61} = \varepsilon_{12} = -1$ and, thus, suspension $\mathcal{S}$ is either convex or concave at edges $e_{61}$ and $e_{12}$ simultaneously, we conclude that

$$\varphi_{61}(x) = \varphi_{12}(x) \qquad \text{for all } x. \tag{4.11}$$

Applying the same arguments to $\mathcal{Q}_1$ treated as being composed of the triangles $\langle \widetilde{N}, \widetilde{p_2}, \widetilde{p_6}\rangle$ and $\langle \widetilde{S}, \widetilde{p_2}, \widetilde{p_6}\rangle$ and taking into account that $\mathcal{Q}_1$ is convex only if it is flat, we get

$$\varphi'_1(x) = -\varphi_1(x) \qquad \text{for all } x. \tag{4.12}$$

Using similar arguments for $\mathcal{Q}_3$ we obtain for all $x$

$$\varphi_{23}(x) = \varphi_{34}(x) \qquad \text{and} \qquad \varphi'_3(x) = -\varphi_3(x). \tag{4.13}$$

$\alpha_j(x)$ **is constant in** $x$ **for** $j = 1, 3, 6,$ **and** $7$ **and is not constant for** $j = 4.$

Substituting (4.9) and (4.11)–(4.13) to (4.1)–(4.7) we obtain

$$\alpha_1(x) = \left(\frac{541419683182996345}{29669606628505029} - \frac{31635727886833754300}{1435086871311616559}\right.$$
$$\left. - \frac{27288800741}{19943076519} + \frac{130585}{9603} - \frac{100}{11} + \frac{310327}{472293}\right)\varphi_{12}(x) + \text{const}, \quad (4.14)$$

$$\alpha_2(x) = -\varphi_3(x) + \varphi_6(x),$$
$$\alpha_3(x) = \text{const},$$
$$\alpha_4(x) = -\varphi_2(x) - \varphi_5(x),$$
$$\alpha_5(x) = \frac{218}{873}(\varphi_1(x) - \varphi_4(x)),$$
$$\alpha_6(x) = \text{const},$$
$$\alpha_7(x) = \text{const}.$$

Note that, due to (3.8), the right-hand side of (4.14) is, in fact, constant in $x$. Hence, $\alpha_j(x)$ is constant in $x$ for $j = 1, 3, 6,$ and $7$. On the other hand we know enough about dihedral angles of $\mathcal{S}$ to prove that $\alpha_4(x)$ or, equivalently, $\varphi_2(x)+\varphi_5(x)$ is not constant.

Recall that Fig. 10 represents suspension $\mathcal{S}$ in a flat position which corresponds to $x = 51$. In the moment, restrict our study by spatial forms of $\mathcal{S}$ which are close enough to that flat position; in particular, assume that

(i) univalent branches of multivalued functions $\varphi_j(x), \varphi'_j(x),$ and $\varphi_{j,j+1}(x)$ are chosen which take values shown on Tab. 6 for $x = 51$ and

(ii) the absolute value of the branch of $\varphi_{12}(x)$ is so small that throughout our discussion there is no necessity to switch to other branches.

Taking into account relations (4.9), (4.11), and (4.13) and using Fig. 10, we draw spherical quadrangles $\mathcal{Q}_2$ and $\mathcal{Q}_5$ on Fig. 16 and 17 respectively.

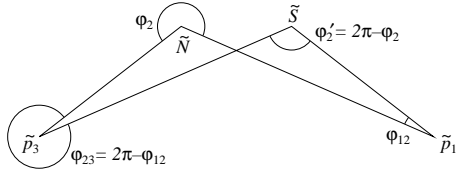 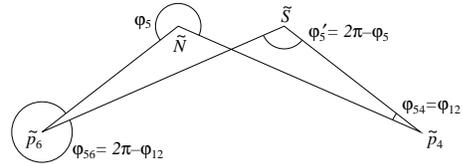

Fig. 16.  Fig. 17.

An obvious consequence of Fig. 16 and 17 is that both $\varphi_2(x)$ and $\varphi_5(x)$ increase as $\varphi_{12}(x)$ increases. Hence, the sum $\varphi_2(x) + \varphi_5(x)$ is not constant. This implies that $\alpha_4(x)$ and Dehn invariant $D_{f_e}$ are not constant and, thus, we have proven the following

**Theorem.** *The suspension $\mathcal{S}$ with a hexagonal equator constructed in § 3 provides us with a counterexample to the Extended Strong Bellows Conjecture.*

**Discussion.**

Note that, for every flexible suspension with a quadrilateral equator (i.e., for every Bricard's flexible octahedron), the Extended Strong Bellows Conjecture holds true [2].

FLEXIBLE SUSPENSIONS WITH A HEXAGONAL EQUATOR 25We conjecture that the Strong Bellows Conjecture is wrong and that a counterexample can be constructed by elimination of self-intersections in a counterexample for the Extended Strong Bellows Conjecture in a way similar to that was used by R. Connelly in [8].

Just note that if the Extended Strong Bellows Conjecture holds true for a polyhedron $\mathcal{P}$ (for example, because $\mathcal{P}$ is not flexible) and is wrong for a polyhedron $\mathcal{Q}$ (for example, consider suspension $\mathcal{S}$ constructed in §3) then it is wrong for a polyhedron $\mathcal{R}$ obtained by glueing $\mathcal{P}$ and $\mathcal{Q}$ along a pair of isometric faces. In particular, if $\mathcal{R}$ is embedded (i.e., has no self-intersections) we arrive to a counterexample to the Strong Bellows Conjecture.

Unfortunately, we cannot realize this idea right now because the suspension $\mathcal{S}$ is more complicated object than the Bricard's flexible octahedron used in [8]. Roughly speaking, Bricard's octahedron in a flat position is a twice-covered polygon, while some parts of $\mathcal{S}$ in a flat position shown on Fig. 9 are four-covered.

At last we wish to underline that flexible suspension $\mathcal{S}$ has surprisingly many hidden symmetries in edges and dihedral angles.

**Acknowledgement.**

The first author expresses his gratitude for the hospitality of the Department of Mathematics of Cornell University during his visit in 2007 when part of this paper was written.## References

[1] R. Alexander, *Lipschitzian mappings and total mean curvature of polyhedral surfaces, I*, Trans. Amer. Math. Soc. **288** (1985), 661–678.
[2] V. Alexandrov, *The Dehn invariants of the Bricard octahedra*, arXiv:0901.2989 [math.MG].
[3] F.J. Almgren and I. Rivin, *The mean curvature integral is invariant under bending*, The Epstein Birthday Schrift, vol. **48**, University of Warwick, 1992, pp. 1–21.
[4] C. Berge, *Graphs and Hypergraphs*, North-Holland, 1973.
[5] M. Berger, *Geometry. Vol. 2*, Springer-Verlag, 1987.
[6] V.G. Boltyanskii, *Hilbert's third problem*, Winston & Sons, 1978.
[7] R. Connelly, *An attack on rigidity I, II*, prepr. (1974), 1–36; 1–28; Russian transl. in the book: A.N. Kolmogorov, S.P. Novikov (eds.), *Issledovaniya po metricheskoj teorii poverkhnostej*, Moscow: Mir, 1980. P. 164–209. Electronic version in English is available at http://www.math.cornell.edu/∼connelly.
[8] R. Connelly, *A counterexample to the rigidity conjecture for polyhedra*, Publ. math. IHES. **47**, 333–338.
[9] R. Connelly, I. Sabitov, and A. Walz, *The bellows conjecture*, Beitr. Algebra Geom. **38** (1997), 1–10.
[10] H. Hadwiger, *Vorlesungen über Inhalt, Oberfläche und Isoperimetrie*, Springer-Verglag, 1957.
[11] F. Harary, *Graph Theory*, Addison Wesley, 1972.
[12] S. Lang, *Elliptic Functions*, Springer-Verlag, 1987.
[13] S.N. Mikhalev, *Some necessary metric conditions for flexibility of suspensions*, Mosc. Univ. Math. Bull. **56** (2001), 14–20.
[14] I.Kh. Sabitov, *Local theory on bendings of surfaces*, Geometry III. Theory of surfaces. Encycl. Math. Sci. **48** (1992), 179–250.
[15] I.Kh. Sabitov, *The volume of a polyhedron as a function of its metric*, Fundam. Prikl. Mat. **2** (1996), 1235–1246. (Russian)
[16] I.Kh. Sabitov, *The volume as a metric invariant of polyhedra*, Discrete Comput. Geom. **20** (1998), 405–425.
[17] R. Walker, *Algebraic Curves*, Dover, 1950.